\documentclass[12pt,letterpaper]{amsart}

\usepackage{amsfonts,amssymb}
\usepackage{amsmath}

\newcommand{\sra}{\rightarrow}             

\renewcommand{\l}{\langle}
\renewcommand{\r}{\rangle}

\newcommand{\N}{\mathbb{N}}  
\newcommand{\Z}{\mathbb{Z}}   
\newcommand{\R}{\mathbb{R}}  
\newcommand{\Q}{\mathbb{Q}}

\newcommand{\mas}[1]{\left\{ #1 \right\}}
\newcommand{\pan}[1]{\left(  #1 \right)}
\newcommand{\hak}[1]{\left[  #1 \right]}

\newcommand{\phib}{\varphi}
\newcommand{\eps}{\varepsilon}

\newcommand{\T}{\mathbb{T}}

\renewcommand{\l}{\langle}
\renewcommand{\r}{\rangle}

\theoremstyle{plain}
\newtheorem{thm}{Theorem}

\newtheorem{lem}[thm]{Lemma}
\newtheorem{pro}[thm]{Proposition}
\newtheorem{cor}[thm]{Corollary}

\newtheorem*{conj}{Conjecture}

\theoremstyle{definition}
\newtheorem{defn}{Definition}

\theoremstyle{remark}
\newtheorem*{rema}{Remark}

\begin{document}
\title[Equidistribution of Polynomial Curves]{Equidistribution of Dilations of Polynomial Curves in Nilmanifolds}
\author{Michael Bj\"orklund and Alexander Fish}
\date{}

\begin{abstract}
In this paper we study the asymptotic behaviour under dilations of probability measures 
supported on polynomial curves in nilmanifolds. We prove, under some mild 
conditions, effective equidistribution of such measures to the Haar measure. 
We also formulate a mean ergodic theorem for $\R^n$-representations on 
Hilbert spaces, restricted to a moving phase of low dimension. 
Furthermore, we bound the necessary dilation of a given smooth curve in $ \R^n$ so that the
canonical projection onto $ \T^n $ is $ \eps$-dense.
\end{abstract}
\maketitle

\section{Introduction}

The concept of mixing in ergodic theory is a fairly well--understood topic from a spectral point of view. In general, mixing of a $\R$-flow $T$ on a 
probability space $(X,\mu)$ asserts that absolutely continuous probability measures $\nu$ converge to $\mu$ under 
the flow, i.e.
\[
(T_t)_*\nu \sra \mu,
\] 
in the weak*--topology on probability measures on $X$. Of course, in general, the attractor of $\mu$ under the dynamics of the flow $T$ can be much larger than
the set of absolutely continuous measures; indeed, it may contain probability measures which are very singular with respect to $\mu$. For an arbitrary mixing $\R$--flow 
there is no precise characterization of the attraction basin of $\mu$. 

In this paper we study the action of dilations on certain classes of singular probability measures on nilmanifolds. This is a slightly different situation than the one
described above, since there is no ergodic $\R$--flow associated to the problem, but the underlying philosophy is the same: Let $ \delta $ denote
a semigroup homomorphism from the additive group $ \R_+$ into the automorphism group of a homogeneous group $G$, which in most of our cases will be $\R^n$.
Let $ \tilde{\mu} $ be a probability measure on $G$, supported on a smooth curve $p$ in $G$, and let $\tilde{\mu}_\lambda$ denote the measure 
$ ( \delta_\lambda)_*\tilde{\mu} $, for $ \lambda > 0 $. Suppose $ \Gamma $ is a cocompact lattice in $G$, and let $ \mu_\lambda $ denote the canonical projection 
of $ \tilde{\mu}_\lambda$ onto the associated manifold. We are interested in conditions on the curve $p$ which assure that the sequence of singular probability 
measures $\mu_\lambda$ converge to the Haar measure, and we want to estimate the speed of convergence. 

In the case of $G = \R^n$, and standard dilations, we give conditions on the curve $p$ which assure fast convergence of $\mu_\lambda$ to the Haar measure 
( in the stronger topology of pseudo--measures, i.e. as functionals on the space of absolutely summable Fourier series ). Furthermore, under some smoothness
conditions on the curve $p$, we bound, for any $ \eps > 0$, the necessary dilation $A_\eps $ from above, so that for all $ \lambda \geq A_\eps$, 
every $\eps$--ball in $\T^n$ is charged by the measures $ \mu_\lambda $.

In the case when $G$ is a Heisenberg group, we prove that the question of equidistribution of a family of orbits \( \{\lambda p\}_{\lambda} \) on a compact quotient $X$ can be reduced
to equidistribution on the horizontal torus of $X$ ( i.e. the quotient of $X$ by the derived group of $G$ ). This result is reminiscent to Leibman's result \cite{Le05} on polynomial sequences in nilmanifolds. 

In section \ref{sec:mean} we prove a general mean ergodic theorem for $\R^n$--actions restricted to a moving phase of dimension one: if $T$ is an ergodic $\R^n$--flow on a countable
generated Borel space $(X,\mu)$, and $p$ is a polynomial curve in $\R^n$ which is not completely contained in an affine hyperplane, then, for every $f \in L^2(X,\mu)$,
\[
\int_{0}^{1} f(T_{\lambda p(t)} x) \, dt \sra \int_{X} f(x) \, d\mu(x),
\]
in the norm of $L^2(X,\mu)$.

The methods used in this paper can also treat analogous equidistribution problems for higher dimensional algebraic varieties.

\textit{Acknowledgment}: 
The authors would like to thank Nimish Shah for mentioning the problem, Manfred Einsiedler for fruitful discussions, and the anonymous referee for suggesting many clarifying remarks on an
earlier version of the paper. 

\section{A Mean Ergodic Theorem} \label{sec:mean}

In this section we will prove a mean ergodic theorem for $\R^n$ actions restricted to a dilated phase of dimension $1$. We will need the following 
simple lemma, which can be found in \cite{Gr}, p. 146.

\begin{lem}[Van der Corput's Lemma]
Suppose $ \phib(t) = \sum_{k = 1}^{d} a_k t^k $ is a polynomial of degree $ d $. Then there is constant $ C $, depending only on $ d $, such that
\[
| \int_{0}^{1} e^{i \phib(t) } \, dt | \leq \frac{C}{\pan{\sum_{k=1}^{d} |a_k|}^{1/d}}.
\]
\end{lem}

For basic facts about Bochner--integration in Hilbert spaces, we refer the reader to \cite{Di} and \cite{Mo01}, 
\begin{thm}
Let $U$ be a unitary $\R^n$--representation, without invariant vectors, on a separable Hilbert space $H$. Let $p$ be a polynomial curve in $\R^n$, defined on the interval 
$ \hak{0,1}$, and not contained in a proper affine subspace of $\R^n$. Define, for $ \lambda > 0 $ and $x \in H$, the Bochner integral
\[
x_\lambda = \int_{0}^{1} U_{\lambda p(t)} x \, dt.
\]
Then, $ x_\lambda \sra 0 $, as $\lambda \sra \infty $, in the norm topology of $H$.
\end{thm}

\begin{proof}
We recall Stone's representation theorem for unitary $\R^n$--repre\-sentations on Hilbert spaces \cite{Sa}: For all $x \in H$, there is a bounded measure $\mu_x$ on $\R^n$ such 
that 
\[
\l x , U_t x \r = \int_{\R^n} e^{i \l y , t \r} \, d\mu_x(y), \quad \forall \: t \in \R^n.
\]
Since $U$ is assumed to be without invariant vectors in $H$, the measure $\mu_x$ does not have an atom at the point $0$. 
Thus, 
\begin{eqnarray}
||x_\lambda||^2 & = & \l \int_{0}^{1} U_{\lambda p(t)} x \, dt , \int_{0}^{1} U_{\lambda p(s)} x \, ds \r  \nonumber \\ 
& = & \int_{0}^{1} \int_{0}^{1} \l U_{\lambda p(t)} x , U_{\lambda p(s)} x \r \, dt ds  \nonumber \\
& = & \int_{0}^{1} \int_{0}^{1} \int_{\R^n} e^{i \lambda \l y , p(s) - p(t) \r} \, d\mu_x(y) dt ds \nonumber \\
& = & \int_{\R^n} \Big|  \int_{0}^{1} e^{i \lambda \l y , p(t) \r } \, dt \Big|^2  \, d\mu_x(y). \nonumber
\end{eqnarray}
Define
\[
\psi_x(\lambda,y) = \Big|  \int_{0}^{1} e^{i \lambda \l y , p(t) \r } \, dt \Big|,
\]
and note that, since $p$ is not contained in a proper affine subspace of $\R^n$, and hence $\l y , p(t) \r $ is a non-trivial polyomial of degree less than or equal to $d$, 
we know that $\psi_x(\lambda,y) \sra 0$, for all $y \neq 0$. Since $ |\psi_x(\lambda,y)| \leq 1$, for all $y \in \R^n$ and $\lambda > 0$, we conclude, by dominated convergence, 
that $ x_\lambda \sra 0 $, in the norm of $H$.
\end{proof}

\begin{rema}
If we assume that the measure $\mu_x$ satisfies an integrability condition of the form, 
\[
\int_{\R^n} \frac{1}{||\phib_y||_1^{2/d}} \, d\mu_x(y) < + \infty,
\]
where $ \phib_y(t) = \l y , p(t) \r $, for $ t \in \hak{0,1}$. Then, if the degree of $p$ is $d$, we have that for all $ \lambda > 0$,
\[
||x_{\lambda}|| \leq \frac{C}{\lambda^{1/d}}, 
\]
for some constant $C$.
\end{rema}
This theorem also extends ( by a straightforward approximation argument ) to isometric representations of $\R^d$ on $L^{p}(X,\mu)$, where $(X,\mu)$ is a 
countably generated Borel space. In this case, it would be interesting to investigate pointwise convergence almost everywhere on 
$(X,\mu)$, and related maximal inequalities.  

\section{Equidistribution of Long Curves on Tori}


In this section we will study equidistribution of dilations of curves in $ \R^n $ projected onto $ \T^n $. Let $ p $ be a smooth curve, defined on the interval
$ \hak{0,1} $. Let $ \mu $ denote the Lebesgue measure on $ \hak{0,1} $, and define $ \mu_{\lambda} $ to be the following measure on $ \T^n $
\[
\int_{\T^n} \phib(x) \, d\mu_\lambda(x) = \int_{0}^{1} \phib(\lambda p(t)) \, dt, \quad \lambda > 0,
\]
where $ \phib $ is a continuous function on $ \T^n = \R^n / \Z^n $. We say that the curves $ \mas{\lambda p}_{\lambda > 0} $ equidistribute in $ \T^n $ if
\[
\int_{\T^n} \phib(x) \, d\mu_\lambda(x) \sra \int_{T^n} \phib(x) \, dm(x), \quad \textrm{as} \:  \lambda \sra \infty
\]
where $ m $ denotes the Haar measure on $ \T^n $. By uniform completeness of exponentials in $ C(\T^n) $ this is equivalent to
\[
\int_{\T^n} e^{2 \pi i \l \nu , x \r } \, d\mu_\lambda(x) \sra 0,
\]
for all $ \nu \in \Z^n \backslash \mas{0} $.

We will restrict our attention to polynomial curves. 

If $ p(t) = \pan{ p_1(t) , \ldots , p_n(t) } $ is a polynomial curve in $ \R^n $, where
\begin{equation}
\label{eq:pol}
p_j(t) = \sum_{k=1}^{d} a_{jk} t^k, \quad j = 1, \ldots , n,
\end{equation}
we let $ A $ denote the $ n \times d $--matrix $ \pan{a_{jk}} $.

We will make effective statements about the convergence of $ \mu_\lambda $ to the Haar measure on $ \T^n $ in terms
of the kernel of  $ A^*$. We let $ || \cdot || $ denote the $ l^1$--norm of the space at hand, and we always assume that $ d \geq n $.

\begin{pro} \label{pro:main}
Suppose the kernel of $ A^* $ is trivial over $ \Q $. Then there is a positive constant $ C $, independent of $ p $, such that for all $ \lambda > 0 $, and $ \nu \in \Z^n \backslash \mas{0} $,
\[
\vert \hat{\mu}_\lambda(\nu)\vert \leq \frac{C}{ \lambda^{1/d} || A^{*} \nu ||^{1/d}}.
\]
\end{pro}
\begin{proof}
This is an immediate consequence of Van der Corput's lemma applied to the function $ \phib(t) = \l \nu , p(t) \r $.
\end{proof}
Note that under the above assumption, no uniformity in $ \nu $ can be asserted. If we assume that the kernel is trivial over $ \R $, more can be said:
\begin{thm}
\label{thm_3}
Suppose the kernel of $ A^* $ is trivial over $ \R $. Then there is a positive constant $ C $, such that for all $ \lambda > 0 $,
\[
\sup_{\nu \neq 0} ||\nu||^{1/d} |\hat{\mu}_{\lambda}(\nu)| \leq \frac{C}{\lambda^{1/d}}.
\]
\end{thm}
In particular, this implies that the measures $ \mu_\lambda $ converge to the Haar measure on $ \T^n $ as pseudo--measures, i.e. in the weak*--topology induced from the space of absolutely summable
Fourier series. By a simple extension \cite{Gr} of van der Corput's lemma, the same theorem holds true if the Lebesgue measure $ \mu $ on $ \hak{0,1} $ is replaced by a compactly supported absolutely continuous
measure on $ \R $ with a sufficiently smooth density. This theorem is similar in spirit to Shah's results \cite{Sh} on equidistribution of smooth curves in the unit tangent bundle of a finite volume hyperbolic manifold
under the action of the geodesic flow.

\begin{proof}
The triviality of the kernel of $ A^* $ over $ \R $ implies the existence of a positive constant $ c  $ such that $ ||A^*x|| \geq c ||x|| $ for all $ x \in \R^n $. The
theorem now follows from proposition \ref{pro:main}.
\end{proof}

We also have the following corollary:

\begin{cor}
\label{thm:dens_equid} The dilations of a polynomial curve  
\( p(t) = (p_1(t),\ldots ,p_n(t)) \) in \( \T^n \) equidistribute if and only
if they become dense.
\end{cor}
\begin{proof}
The only direction which requires a proof is ``\( \Leftarrow\)''.  
If the dilations do not equidistribute, then, by proposition \ref{pro:main}, there must be a linear dependence 
between the coefficients of the polynomials \(p_1,\ldots, p_n\). This implies that there exists 
\( m \in \Z^n \setminus \{0\} \) such that for all $ \lambda \in \R $ and $ t \in \hak{0,1} $, 
\[
\l m,\lambda p(t) \r = \lambda \l A^*m , \gamma(t) \r = 0,
\]
where $ \gamma(t) = ( t , t^2 , \ldots , t^d ) $.

This means that all the curves \( \lambda p \) lie in the subtorus 
\[
L_m = \{ x \in \T^n \, | \, \l m,x\r = 0 \}.
\] 
Thus, in this case, dilations do not spread out densely in $\T^n$.
\end{proof}

\section{Equidistribution on Heisenberg Nilmanifolds}
\label{sec3}

The Heisenberg groups $H_n(R)$, where $ R $ is a unital commutative ring, are defined as follows:
\[
H_n(R) =
\left\{
\left(
\begin{array}{ccccc}
1 & x_1 & \ldots & x_n & z \\
 & \ddots & & & y_1 \\
 & & \ddots & 0 & \vdots \\
 & 0 & & \ddots & y_n \\
& & & & 1
\end{array}
\right)
\: | \: x,y \in R^n, z \in R
\right\},
\]
It is straightforward to prove that $ H_n(R) $ is a two--step nilpotent group. If $ R = \R $, $ H_n(\R) $ is a Lie group, with the smooth structure induced from the Lie algebra (see \cite{Ta}, p. 42),
and $ H_n(\Z) $ is a cocompact lattice in $ H_n(\R)$. We let $ X_n = H_n(\R) / H_n(\Z) $ denote the associated nilmanifold, and we let $ m $ denote the unique $ H_n(\R)$--invariant measure on $ X_n $.
The Laplace operator on $ X_n $ is given by
\[
\Delta = \sum_{j=1}^{n} ( D_j^2 + D_j'^2) + \partial_z^2,
\]
where $ D_j  = \partial_{x_j} $ and $ D'_j = \partial_{y_j} + x_j \partial_z $.

The eigenfunctions and eigenvalues of $ \Delta $ on $ X_n $ were determined by Deninger and Singhof in \cite{DeSi84}. We briefly recall their result: The spectrum of $ \Delta $ decomposes
into two parts; the first part is parameterized by $ k , h \in \Z^n $, with eigenfunctions
\[
f_{k,h}(x,y,z) = e^{2 \pi i( \l k , x \r + \l h , y \r )}
\]
and eigenvalues
\[
\lambda_{k,h} = -4\pi^2(||k||^2 + ||h||^2).
\]
The second part of the spectrum consists of eigenfunctions of the form,
\[
g_{q,m,h}(x,y,z) = e^{2 \pi i(mz + \l q , y \r)}
\prod_{j=1}^{n} \hak{ \sum_{k \in \Z} F_{h_j}\pan{ \sqrt{2\pi |m|} \pan{x_j + \frac{q_j}{m}+k}} e^{2\pi i k m y_j}}
\]
where $ m \in \Z \backslash \mas{0}, q \in \Z^n $ and $ h = ( h_1 , \ldots , h_n ) \in \N_0^n $, with eigenvalues
\[
\lambda_{q,m,h} = -2\pi|m|(2h_1 + \ldots + 2h_n + n + 2\pi |m|).
\]
Here, $ \mas{F_\nu}_{\nu \geq 0} $ denotes the Hermite functions,
\[
F_\nu(t) = (-1)^{\nu} e^{t^2/2} \frac{d^\nu}{dt^\nu}e^{-t^2}, \quad \nu \geq 0.
\]
Note that the $ H_n $--invariant measure $ m $ on $ X_n $ is completely determined by the equations
\[
\int_{X_n} f_{k,h} \, dm = 0 \quad \textrm{and} \quad \int_{X_n} g_{q,m,h} \, dm = 0.
\]
Indeed, the linear span of eigenfunctions of $ \Delta $ is dense in the space of continuous functions on $ X_n $, see  \cite{Wa}, p. 256. Thus, a necessary and 
sufficient condition for a family of measures $ \mu_\lambda $ to converge to $ m $ in the weak*--topology of the space of probability measures on $ X_n $, is
\[
\int_{X_n} f_{k,h} \, d\mu_\lambda \sra 0 \quad \textrm{and} \quad \int_{X_n} g_{q,m,h'} \, d\mu_\lambda \sra 0,
\]
for all $ (k,h) \in \Z^{2n} \backslash \mas{0} $, $ q \in \Z^{n}  $, $ m \in \Z \backslash \mas{0} $ and $ h' \in \N_0^{n}  $.

The Heisenberg groups admit natural families of dilations, similar to the Euclidean case. We define, for $ \lambda > 0 $,
\[
\lambda(x,y,z) = ( \lambda^{1/2} x , \lambda^{1/2} y , \lambda z ).
\]
Let
\[
p(t) = ( a(t) , b(t) , c(t) ), \quad t \in \hak{0,1},
\]
where $ a $ and $ b $ are polynomial curves in $ \R^n $, and $ c $ is a polynomial.
 As before, we let $ \mu_{\lambda} $
denote the push--forward of the Lebesgue measure on the interval $ \hak{0,1} $ under the composite of $ \lambda p $ into $ H_n $ and the canonical projection
of $ H_n $ onto the nilmanifold $ X_n $. The horizontal torus of $ X_n $ is the torus $ \T^{2n} $ generated by the $ x $ and $ y $--coordinates modulo $ \Z^{2n} $.
The following theorem should be thought of as an analogue of Leibman's theorem \cite{Le05}:

\begin{thm}
\label{thm:main_nilman} Suppose $ p $ is a polynomial curve in $
H_n(\R) $. Then, 
the dilations of $ p $ equidistribute in the nilmanifold 
$ X_n = H_n(\R) / H_n(\Z) $ if and only if they equidistribute in
the horizontal torus of $ X_n $.
\end{thm}

\begin{proof}
We will first prove that
\[
I(\lambda) = \int_{X_n} g_{q,m,h'} \, d\mu_\lambda \sra 0, \quad \lambda \sra \infty,
\]
for all $ q \in \Z^{n} 
$, $ m \in \Z \backslash \mas{0} $ and $ h' \in \N_0^{n}  $.

We define $G(x,y) = G_1(x_1,y_1) \cdots G_n(x_n,y_n)$, and
\[
G_j(x_j,y_j) = \sum_{k \in \Z} F_{h_j} \Big( \sqrt{2\pi |m|} \Big( x_j + \frac{q_j}{m} + k \Big) \Big) e^{2\pi i (km+q_j)y_j}.
\]
The general case can be deduced from the case when \( c \) is non--negative and invertible on \( [0,1]\) with \( c(0) = 0 \). Let \( [0,\beta] = c ([0,1]) \). 
In this new notation we are interested in the asymptotics of the integral
\[
 I(\lambda) = \int_{0}^{1} e^{2 \pi i m \lambda c(t)} G(\lambda^{1/2}a(t),\lambda^{1/2}b(t))dt,
\]
when \( \lambda \to \infty \). We make the variable substitution \( u = c(t) \):
\[
 I(\lambda) = \int_{0}^{\beta} \frac{e^{2 \pi i m \lambda u} G(\lambda^{1/2}a(c^{-1}(u)),\lambda^{1/2}b(c^{-1}(u)))}{c'(c^{-1}(u))}du.
\]
Note that the function 
\[
 H_{\lambda}(u) = G(\lambda^{1/2}a(c^{-1}(u)),\lambda^{1/2}b(c^{-1}(u)))
\]
is bounded as a function in \( \lambda \) and \( u \), and if we let \( \omega(u) = c'(c^{-1}(u)) \), we have 
\[
 \int_{0}^{\beta} \frac{1}{|\omega(u)|}du < +\infty.
\]
Thus, for any \( \varepsilon > 0 \) there is \( \delta > 0 \) such that 
\[
 \left| \int_{0}^{\delta} \frac{1}{\omega(u)}\right| < \varepsilon.
\]
Hence, since \( H_{\lambda}(u) \) as a function of \( \lambda \) is bounded for \( u \in [0,\delta] \), to prove that \( I(\lambda) \) converges to zero as \( \lambda \to +\infty \), 
it suffices to estimate the following integral:
\begin{eqnarray*}
&& \int_{\delta}^{\beta} \frac{e^{2 \pi i m \lambda u} G(\lambda^{1/2}a(c^{-1}(u)),\lambda^{1/2}b(c^{-1}(u)))}{c'(c^{-1}(u))}du \\
&=& \frac{1}{2 \pi i m \lambda} \int_{\delta}^{\beta} \frac{(e^{2 \pi i m \lambda u})' G(\lambda^{1/2}a(c^{-1}(u)),\lambda^{1/2}b(c^{-1}(u)))}{c'(c^{-1}(u))}du \\
&=& \frac{o(\lambda)}{\lambda} - \frac{1}{2 \pi i m \lambda} \int_{\delta}^{\beta} e^{2 \pi i m \lambda u}\left(\frac{H_{\lambda}'(u)}{\omega(u)} - H_{\lambda}(u) 
\frac{\omega'(u)}{\omega(u)^2} \right) du.
\end{eqnarray*}
where $'$ denotes differentiation with respect to the variable $u$. 

The last term is clearly \( o(\lambda) \). Since \( H_{\lambda}'(u) = \lambda^{1/2} K(\lambda,u) \) for some function $K$
which is bounded in \( u \) and \( \lambda \), we conclude that 
\[
 \Big| \int_{\delta}^{1} e^{2\pi i m \lambda u} \frac{H_\lambda'(u)}{\omega(u)} \Big| = O(\lambda^{1/2}),
\]
which proves that $ \limsup_{\lambda \sra \infty} |I(\lambda)| < \eps $. Since $\eps$ is arbitrary, we are done.
 
On the other hand,
\[
\int_{X_n} f_{k,h} \, d\mu_\lambda \sra 0, \quad \lambda \sra \infty,
\]
is equivalent to the statement that the projections  of $ \lambda p $ onto $ X_n $ equidistribute on the horizontal torus.
\end{proof}

As a corollary of Theorems \ref{thm:main_nilman} and
corollary \ref{thm:dens_equid} we get the following result.

\begin{thm}
Suppose $ p $ is a polynomial curve in $ H_n(\R) $. The dilations of
$ p $ equidistribute in the nilmanifold $ X_n = H_n(\R) / H_n(\Z) $
if and only if they become dense in the horizontal torus of $ X_n $.
\end{thm}

 For a general nilmanifold $ X = N / \Gamma $, where $ N $ is a
nilpotent Lie group, and $ \Gamma $ a cocompact discrete subgroup of
$ N $, the horizontal torus is defined to be the manifold $ N /
\hak{N,N}\Gamma $. We conjecture that the above phenomena is not
restricted to Heisenberg manifolds, but occurs in any nilmanifold
which admits dilations.

\begin{conj}
Suppose $ X $ is a compact nilmanifold, covered by a nilpotent Lie group $N$ which admits dilations. If $p$ is a polynomial curve in $N$, the canonical projections onto $X$ of the dilations of $p$
equidistribute, if and only if they equidistribute in the horizontal torus of $X$.
\end{conj}

This conjecture is similar to Leibman's theorem \cite{Le05} on equidistribution of polynomial sequences in nilmanifolds. 

\section{Density of Polynomial Curves}

In this section we will address the following question: \emph{Given a polynomial curve $ p $ in $ \R^n $, defined on $ \hak{0,1}$, and $ \eps > 0 $; find
a number  $ A(\eps) $ such that for all $ \lambda \geq A(\eps) $, the canonical projection of the orbit $ \lambda p(\hak{0,1}) $ onto $ \T^n $ is $ \eps $--dense}.
The answer will of course depend on certain diophantine properties of the coefficients of the polynomial curve. We define, for $ z \in \T $,
\[
\Vert z \Vert_{\T} = \inf_{n \in \Z} | z - n |.
\]
We will need the following definition:
\begin{defn}
Let \( c >0 \) and \( q \in \N \). A vector $ x $ in $ \R^n $ is called \emph{$ ( c , q ) $--badly approximable by rationals} \mbox{( BAP )} if  for every
$ \nu $ in $ \Z^n \backslash \mas{0} $ we have
\[
\Vert \l x , \nu \r \Vert_{\T} > \frac{c}{\Vert \nu \Vert^{q}}.
\]
\end{defn}
It follows from the pigheonhole principle that, for every dimension \( n \in \N \), there exists a constant \(c_n > 0 \) such that every vector  \( x \in \R^n \) is not \( (c_n,n) \)-BAP. For example \( c_1 = 1 \).  In general, \( c_n = n^{\frac{n}2}\).
We say that a vector \( x \in \R^n \) is \textit{badly approximable} if there exists \( c > 0 \) such that \( x \) is \((c,n)\)-badly approximable. Denote by \( \mathcal{B}A \) the set of all badly approximable vectors. 

The well known Khintchine's theorem states that the set \( \mathcal{B}A \) has zero Lebesgue measure. Nevertheless, it is a quite large set. Schmidt's theorem \cite{Sc} states that \( \mathcal{B}A \) is a \textit{thick} set, i.e., for every open set \( W \subset \R^n \) the Hausdorff dimension of \( \mathcal{B}A \cap W \) is equal to \( n \).
One more result should be mentioned in the context of a largeness of \( (c,q)\)-BAP sets. For every \( n \in \N \) and every \( \eps > 0 \) the set of vectors in \( \R^n \) which are \((c,n+\eps) \)-BAP for some \( c > 0 \) (which might depend on \( x \)) has full Lebesgue measure.


We will prove:
\begin{thm} \label{main:eps}
Suppose $f$ is a curve in $ \R^n $ defined on $\hak{0,1}$ such that \( f(0) = 0 \) and \( f\) is twice continuously differentiable at \( 0 \). Assume that \( f'(0) \) is \( (c,q) \)-badly approximable by rationals. Then there exists \( \eps_0 > 0 \) such that for every \( \eps > 0 \) there exists a constant $ C $, depending only on $ c, q, n, f $, such that for every $ \lambda \geq A(\eps) $, where
\[
A(\eps) = C\eps^{-\left((q+1)n(n+1) + q \sum_{i=0}^{n-2} (i+1)(2n-i) + 1\right)},
\]
the projected orbit $ \lambda f(\hak{0,1}) $ onto $ \T^n $ is $ \eps $--dense.
\end{thm}
\textbf{Remark:}
Note that in the most interesting case \( q = n \) there exists \( \alpha > 0 \) (which does not depend on \( f \)) such that 
\[
 A(\eps) = C \eps^{- \alpha n^4 }.
\]

As a corollary of Theorem \ref{main:eps} we get an upper bound on the dilation \(A(\eps) \) which will suffice for \(\eps\)-density of a dilated polynomial curve \( p(t) \) under some diophantine conditions on the coefficients.
\begin{cor}
Let \( p(t) = (p_1(t),\ldots,p_n(t)) \) be a polynomial curve in \( \R^n \) given by equations (\ref{eq:pol}).
Assume that the vector \( (a_{11},a_{21},\ldots,a_{n1}) \) is \((c,q)\)-badly approximable by rationals. Then 
there exists \( C > 0 \) (which depends on \( p, n, c, q\)) such that for every \( \lambda \geq A(\eps) \), where
\[
A(\eps) = \frac{C}{\eps^{2 \left((q+1)\frac{n(n+1)}{2} + q \sum_{i=0}^{n-2} \frac{(i+1)(2n-i)}{2}\right) + 1}},
\]
the projected orbit \( \lambda p([0,1]) \) onto \(\T^n \) is \( \eps \)-dense. 
\end{cor}

Let $x \in \R^n $. Define
\[
T_{x}(\eps) = \inf\mas{ N \geq 1 \: | \: \textrm{the orbit} \: \mas{mx}_{m=1}^N \,\,  \textrm{is $\eps$--dense in $ \T^n $ } }.
\]

Theorem \ref{main:eps} follows immediately form the following upper bound on  \( T_{x}(\eps) \), when $x$ is BAP.
\begin{lem}
\label{kronecker_type_lemma}
Assume \( {a} = (a_1,\ldots,a_n) \in \R^n\) is \((c,q)\)-badly approximable by rationals. Then for \( n > 1 \) we have the following estimate: 
\[ T_{{a}}(\eps) <<  \left( \frac{1}{\eps}\right)^{(q+1)\frac{n(n+1)}{2} + q \sum_{i=0}^{n-2} \frac{(i+1)(2n-i)}{2}} =: \Lambda(\eps), 
\]
for \( n = 1 \) we have the estimate:
\[
 T_{{a}}(\eps) \leq \frac{1}{c \eps^{q+1}}.
\]

\vspace{0.1in}

\noindent
\textnormal{The expression  ``\( e_1(x) << e_2(x) \)'' means that there exists a constant \( C > 0 \) (which does not depend on \( x \)) such that \( e_1(x) \leq C e_2(x), \forall x > 0 \).}
\end{lem}
\begin{proof}[Proof of Theorem \ref{main:eps}.]
Let  \( f_{\lambda}(t) =  \lambda f(t) \). We will find an upper bound on \( A(2\eps) \) as a function of \( \Lambda(\eps) \) obtained in Lemma \ref{kronecker_type_lemma}.

Let \( \lambda > 0 \) an let \( a = f'(0) \). By the assumption on \( f(t) \) and Lemma \ref{kronecker_type_lemma} we know that there exists \( D > 0\) such that for sufficiently small  $\eps_0 > 0 $, 
\begin{enumerate}
 \item \( \sup_{t \leq D \eps_0^{\frac{1}2}} \Vert f_{\lambda}(t) - a \lambda t \Vert \leq \lambda \eps_0.\)
\item  \(\{ \lambda a t \, | \, t \leq D \eps_0^{1/2} \} \) is \( \eps \)-dense in \( \T^n \) provided that \( \lambda D \eps_0^{1/2} \geq \Lambda(\eps) \).  
\end{enumerate}
Thus, any \( \lambda \) which satisfies
\begin{enumerate}
 \item \( \lambda \eps_0 \leq \eps \).
\item \label{eq_second} \( \lambda D \eps_0^{1/2} \geq \Lambda(\eps) \).
\end{enumerate}
gives an upper bound on \( A(2\eps) \).
We are interested in a minimal possible \( \lambda \) which satisfies the properties listed above, so we assume that 
\[
 \lambda \eps_0 = \eps.
\]
We deduce \( \eps_0 = \frac{\eps}{\lambda} \). We substitute the obtained expression for \( \eps_0 \) into the equation (\ref{eq_second}):
\[
 \lambda^{1/2} \eps^{1/2} D \geq \Lambda(\eps).
\]
 The latter inequality implies that 
\[
 \lambda \geq D^{-2}\frac{\Lambda(\eps)^{2}}{\eps}.
\]
The latter implies the conclusion of the theorem.
\end{proof}

\begin{proof}[Proof of Lemma \ref{kronecker_type_lemma}.]
First, we prove the case \( n = 1 \).

It is clear by the pigeonhole principle that there exists \( m \in \N, m \leq \frac{1}{\eps} \) such that \( \Vert ma \Vert_{\T} \leq \eps \). We know that \( a \) is \( (c,q)\)-BAP, thus we have
\[
 \Vert ma \Vert_{\T} \geq \frac{c}{m^q} \geq c \eps^q.
\]
Therefore the first \( \frac{1}{\eps} \frac{1}{c\eps^q} \) iterations of \( a \) will form an \( \eps\)-net in \( \T \).
 
\textbf{The case \( n > 1 \)}.

If we divide \( \T^n \) into \( c_1 \frac{1}{\varepsilon^n} \) boxes then  by pigeonhole principle we can find \( m \leq  N = c_1 \frac{1}{\varepsilon^n} \) such that $ \Vert m {a} \Vert_{\T^n} \leq \varepsilon \). Notice that \( c_1 \) depends only on \( n \).

Let $ {b} = (b_1,\ldots,b_n) $ denote the vector  $ m {a} $. Note that if we can show that the finite orbit \( ({b}, 2{b},\ldots,M{b}) \) is \( \varepsilon \)-dense in \( \T^n \) for some \( M \in \N \), then it follows that
\[
T_{{a}}(\varepsilon) \leq c_1 \frac{M}{\varepsilon^n}.
\]
The new vector \( {b} \) is (\(c/m^q ,q\))-BAP, since
\[
 \Vert \langle ma, \nu \rangle \Vert_{\T} = \Vert \langle a, m \nu \rangle \Vert_{\T} \geq \frac{c}{{(m \Vert \nu \Vert)}^q}.
\]
 Because of the estimate on \( m \) we conclude that \( b \) is \((\frac{c\eps^{nq}}{c_1^q},q)\)-BAP.

Therefore \( b_i \neq 0 \) for all \( i \in \{1,2,\ldots,n\} \). Permutation of the coordinates in the vector \( {b} \) does not change \( T_{{b}}(\varepsilon)\), thus without loss of generality we can assume that \( b_1 \geq b_i, \forall \: 2 \leq i \leq n \).


The trajectory of the flow $ t \mapsto t {b} $ in $  \T^n $ hits the $ (n-1) $-dimensional torus
\[
\T_0 = \{(0,x_2,\ldots,x_n)\in \T^n\}
\]
for $ t \in \frac{1}{b_1} \N $.

We will find an upper bound \( U \) on the number of consecutive visits of \( \T_0 \) by the flow which guarantees that the hit points in \( \T_0 \) constitute  an \( \eps\)-net.

From (\(c/m^q ,q\))-BAP of the vector \( {b} \) we can conclude that \( b_1 \geq c/m^q \). We use that 
\( m \leq c_1 / \eps^n \) and we get \( b_1 \geq \frac{c \eps^{nq}}{c_1^q}\).

Let us prove that the following upper bound 
\[
 T_a(3 \eps) \leq \frac{c_1 U}{{\eps}^n b_1} \leq \frac{c_1^{q+1} U}{c {\eps}^{n(q+1)}}
\]
holds.
First, for any point \( z \in \T_0 \) there exists a point \( y = (k/b_1) b \in \T_0 \) with \( k \in \N, k \leq U \) such that \( dist_{\T_0}(z,y) < \eps\). 
Take \( \ell \in \N \) be equal to \( \ell = \lfloor{k/b_1}\rfloor \).
So \( \ell \leq \frac{U}{b_1} \). Then  
\[ dist_{\T^n}(\ell b, z)  \leq dist_{\T^n}(\ell b,y) + dist_{\T^n}(y,z) \leq 2 \eps. \]
  Let \( z' = (z_1,z_2,\ldots,z_n) \in \T^n \) be an arbitrary point. 
If we look at the finite part of the orbit of \( \{ tb \, | \, t \leq U \} \) intersected with the affine \((n-1)\)-dimensional torus \(\T_{z_1} = \{ (z_1,x_2,\ldots,x_n) \in \T^n \} \) then we will see a \( 2\eps \)-net of
 \( \T_{z_1} \). We move the picture
at the zero level by an isometry to the level of \( z_1 \):
\[ \T_0 \cap \{bt \, | \, t \leq U \} + \frac{z_1}b_1 b \subset \T_{z_1} \cap \{bt \, | \, t \leq U + z_1/b_1\}.   \]
 We might  miss one point from the zero level. A possible missed point contributes by adding one more \( \eps \). By repeating the previous considerations that have been used at level zero we conclude that we can find \( \ell \in \N, \ell \leq \frac{U}{b_1} \) such that 
\[
 dist_{\T^n}(\ell b, z') \leq 3 \eps.
\]

At the moment, we give an argument for finding an upper bound on \( U \). 
Look at the new vector \( {w} = ( \frac{b_2}{b_1},\ldots,\frac{b_n}{b_1}) \in \T^{n-1} \).
Notice that \( (b_2,\ldots,b_n) \) is \((c\eps^{nq} / (c_1^q),q)\)-BAP vector (because \( b \) has the same parameters):
\[
 \Vert \langle (b_2,\ldots,b_n), (\nu_2,\ldots,\nu_n) \rangle \Vert_{\T} = \Vert \langle b, (0,\nu_2,\ldots,\nu_n) \rangle \Vert_{\T} \geq 
\]
\[
\frac{c\eps^{nq} / c_1^q}{\Vert (0,\nu_2,\ldots,\nu_n)\Vert^q} = \frac{c\eps^{nq} / c_1^q}{\Vert (\nu_2,\ldots,\nu_n)\
\Vert^q}.
\]

We claim that \( w \) is \(\left(\frac{c \eps^{nq}}{c_1^{q}(n+1)^q},q\right)\)-BAP.
Let \( \nu \in \Z^{n-1} \setminus\{0\} \). 
\[
 \Vert \langle w,\nu \rangle \Vert_{\T} = 
\min_{m \in \mathbb{Z}} \left| \langle w, \nu \rangle - m \right| = 
\frac{1}{b_1} \min_{m \in \Z} \left| \langle (b_2,\ldots,b_n), \nu \rangle - b_1 m \right| =
\]
\[
\frac{1}{b_1} \min_{|m| \leq |\langle w, \nu\rangle|  +1} \left| \langle b,(m,\nu)\rangle\right| \geq 
 \frac{1}{b_1}\min_{|m| \leq |\langle w, \nu\rangle|  +1} \frac{c \eps^{nq}}{c_1^q \Vert (m,\nu)\Vert^q} = I.
\]
We used the notation \( (m,\nu) \) for a vector in \( \Z^n \) with the first coordinate equal to \( m \) and the rest of coordinates are equal to the vector \( \nu \).
We use the following estimate to bound \( m \):
\[
 |m| \leq \left| \langle w,\nu\rangle\right| +1 \leq \Vert w \Vert \Vert \nu \Vert +1 \leq \sqrt{n-1} \Vert \nu \Vert + 1 \leq n \Vert \nu \Vert.
\]
Thus
\[
 \Vert (m,\nu) \Vert \leq |m| + \Vert \nu \Vert \leq (n+1) \Vert \nu \Vert.
\]
Finally.
\[
I \geq \frac{1}{b_1} \frac{c \eps^{nq}}{c_1^q (n+1)^q \Vert \nu \Vert^q} \geq
\frac{c \eps^{nq}}{c_1^q (n+1)^q \Vert \nu \Vert^q}.
\]

By the definition of \( T_{\omega}(\eps) \) we get the bound:
\[
U \leq T_{{w}}(\varepsilon).
\]
Therefore, we have a recursive formula:
\begin{equation}
\label{recursion}
 T_{{a}}(3\varepsilon)  \leq \frac{c_1^{q+1}}{c \eps^{n(q+1)}} T_{{w}}(\varepsilon),
\end{equation}
where the new vector \( {w} \) is \(\left(\frac{c \eps^{nq}}{c_1^{q}(n+1)^q},q\right)\)-BAP.
\vspace{0.2in}

We have already shown (case \( n =1 \)):
\textit{If \( \alpha \in \T \) is \( (c,q) \)-BAP then
\begin{equation}
\label{torus}
T_{\alpha}(\varepsilon) \leq \frac{1}{c\varepsilon^{q+1}}.
\end{equation}}
\vspace{0.2in}

Let \( f(\eps) \) be an arbitrary positive function of \( \eps \).
We  prove by induction that if \( a \in \R^n \) is \( (f(\eps),q) \)-BAP vector then the following estimate holds:
\[
 T_{{a}}(\varepsilon) << \left( \frac{1}{f(\eps)}\right)^n \left( \frac{1}{\eps} \right)^{(q+1)\frac{n(n+1)}{2} + q \sum_{i=0}^{n-2} \frac{(i+1)(2n-i)}{2}  },
\]
where for \( n  = 1 \) we define \( \sum_{i=0}^{n-2} \frac{(i+1)(2n-i)}{2} = 0\). 
For \( n = 1 \) the formula coincides with formula (\ref{torus}).
For \( n+1 \), for \( a \in \R^{n+1} \) which is \((f(\eps),q)\)-BAP we have by formula (\ref{recursion}):
\begin{equation}
\label{eq_1}
 T_{{a}}(\varepsilon) << \frac{1}{f(\eps) \eps^{(n+1)(q+1)}} T_{{w}}(\varepsilon),
\end{equation}
where \( w \in \R^n \) is \( \left(\frac{f(\eps) \eps^{(n+1)q}}{c_2^{q}(n+2)^q},q\right) \)-BAP (\(c_2 \) depends only on \( n+1 \)). By the induction hypothesis we conclude:
\begin{equation}
\label{eq_2}
 T_{w}(\eps) << \left( \frac{1}{f(\eps)}\right)^n \left( \frac{1}{\eps}\right)^{(n+1)nq} \left( \frac{1}{\eps} \right)^{(q+1)\frac{n(n+1)}{2} + q \sum_{i=0}^{n-2} \frac{(i+1)(2n-i)}{2}  }
\end{equation}
We merge equations (\ref{eq_1}) and (\ref{eq_2}) together and get
\[
 T_{{a}}(\varepsilon) << \left( \frac{1}{f(\eps)}\right)^{n+1} \left( \frac{1}{\eps}\right)^{(n+1)nq + (q+1)\frac{n(n+1)}{2} + q \sum_{i=0}^{n-2} \frac{(i+1)(2n-i)}{2} + (n+1)(q+1)} =
\]
\[
 \left( \frac{1}{f(\eps)}\right)^{n+1} \left( \frac{1}{\eps}\right)^{(q+1)\frac{(n+1)(n+2)}{2} + q \sum_{i=0}^{n-1} \frac{(i+1)(2(n+1)-i)}{2}}.
\]

Therefore, if \( a \in \R^n \) is a \( (c,q) \)-BAP vector then 
\[
 T_a(\eps) << \left( \frac{1}{\eps}\right)^{(q+1)\frac{n(n+1)}{2} + q \sum_{i=0}^{n-2} \frac{(i+1)(2n-i)}{2} }.
\]
\end{proof}

\end{document}